\documentclass[12pt,reqno]{amsart}

\usepackage{pb-diagram} 
\usepackage{amsmath,amssymb}
\usepackage{amsthm}
\usepackage[colorlinks = true,
            linkcolor = red,
            urlcolor  = blue,
            citecolor = blue,
            anchorcolor = blue]{hyperref}
\usepackage[active]{srcltx}

\usepackage{enumitem}   

\newtheorem{theorem}{Theorem}

\newfont{\bb}{msbm10 at 12pt}

\def\SB{\mathbf{S}\Omega}

\def\<{\langle}     
\def\>{\rangle}

\def\Tr{{\rm Tr}}

\newcommand{\bal}{\begin{align}}      \newcommand{\eal}{\end{align}}
\newcommand{\ba}{\begin{array}}      \newcommand{\ea}{\end{array}}
\newcommand{\bc}{\begin{center}}     \newcommand{\ec}{\end{center}}
\newcommand{\be}{\begin{enumerate}}  \newcommand{\ee}{\end{enumerate}}
\newcommand{\beQ}{\begin{eqnarray*}} \newcommand{\eeQ}{\end{eqnarray*}}
\newcommand{\bi}{\begin{itemize}}    \newcommand{\ei}{\end{itemize}}
\newcommand{\bt}{\begin{tabular}}    \newcommand{\et}{\end{tabular}}
\newcommand{\bdm}{\begin{displaymath}} \newcommand{\edm}{\end{displaymath}}

\newcommand{\D}{D\!\!\!\!/\,}
\newcommand{\nb}{\nabla\!\!\!\!/\,}
\newcommand{\RSB}{S\!\!\!\!/\,}
\newcommand{\mult}{c\!\!\!/}

\begin{document}

\title[Nonexistence of DEC spin fill-ins]{Nonexistence of DEC spin fill-ins}       

\author{Simon Raulot}
\address[Simon Raulot]{Laboratoire de Math\'ematiques R. Salem
UMR $6085$ CNRS-Universit\'e de Rouen
Avenue de l'Universit\'e, BP.$12$
Technop\^ole du Madrillet
$76801$ Saint-\'Etienne-du-Rouvray, France.}
\email{simon.raulot@univ-rouen.fr}

\begin{abstract}
In this note, we show that a closed spin Riemannian manifold does not admit a spin fill-in satisfying the dominant energy condition (DEC) if a certain generalized mean curvature function is point-wise large. 
\end{abstract}

\date{\today}   
\maketitle
\pagenumbering{arabic}

Given a closed $n$-dimensional Riemannian manifold $(\Sigma,\gamma)$, it is a very interesting question to know whether there exists an $(n+1)$-dimensional compact Riemannian manifold $(\Omega,g)$ with nonnegative scalar curvature (NNSC) whose boundary is isometric to $(\Sigma,\gamma)$. If so, the Riemannian manifold is called a fill-in of $(\Sigma,\gamma)$ and the set of all such fill-ins, referred to as NNSC fill-ins of $(\Sigma,\gamma)$, is denoted by $\mathcal{F}(\Sigma,\gamma)$. The problem of the existence of such fill-ins has recently been solved by Shi, Wang and Wei \cite{ShiWangWei1} where it is shown that if $\Sigma$ is the boundary of an $(n+1)$-dimensional compact manifold $\Omega$ then, any metric $\gamma$ on $\Sigma$ can be extended to a Riemannian metric $g$ on $\Omega$ with positive scalar curvature.

One can also try to find a fill-in whose mean curvature is prescribed by a smooth function $H$ on $\Sigma$. This problem is tightly related to the Bartnik quasi-local mass \cite{Bartnik2} and a triplet $(\Sigma,\gamma,H)$ is then usually called a Bartnik data. 

In \cite{Miao3}, Miao proved that if $\Sigma$ is the boundary of some compact $(n+1)$-dimensional manifold $\Omega$, then given any Riemannian metric $\gamma$ on $\Sigma$, there exists a constant $H_0$, depending on $\gamma$ and $\Omega$, such that, if $\min_\Sigma H\geq H_0$, there does not exist NNSC fill-ins of $(\Sigma,\gamma,H)$. The proof makes use of the work of Shi, Wang and Wei \cite{ShiWangWei1} and of Schoen and Yau's results on closed manifolds \cite{SchoenYau7,SchoenYau3}. Such a result is also obtained for fill-ins with a negative scalar curvature lower bound. 

This fact was previously demonstrated by Gromov \cite{Gromov1} for spin manifolds. More precisely, he showed that if $(\Omega,g)$ is a NNSC spin fill-in of $(\Sigma,\gamma,H)$, then  
\begin{eqnarray*}\label{GromovInequality}
\min_\Sigma H\leq \frac{n}{{\rm Rad}(\Sigma,\gamma)}
\end{eqnarray*}
where ${\rm Rad}(\Sigma,\gamma)$ is a constant only depending on $(\Sigma,\gamma)$ and known as the hyperspherical radius of $(\Sigma,\gamma)$. The main remark here, which is the key point of our argument, is that the nonexistence of NNSC spin fill-ins can be obtained by taking another route. Indeed, it can be deduced from from an eigenvalue estimate on the first eigenvalue $\lambda_1(\Sigma,\gamma)$ of the Dirac operator of $(\Sigma,\gamma)$ proved by Hijazi, Montiel and Zhang \cite{HijaziMontielZhang1} and which states that 
\begin{eqnarray}\label{HMZInequality}
\min_\Sigma H\leq 2|\lambda_1(\Sigma,\gamma)|.
\end{eqnarray} 
Since $\lambda_1(\Sigma,\gamma)$ depends only on $\Sigma$, $\gamma$ and the involved spin structure, the nonexistence of NNSC spin fill-ins is a direct consequence of (\ref{HMZInequality}). This inequality is sharp since round balls in the Euclidean space satisfy the equality case of (\ref{HMZInequality}). Note that a similar result can be deduced for spin fill-ins with a negative scalar curvature lower bound from \cite{HijaziMontielRoldan} as in \cite{Gromov1}. 

In this note, we use the aforementioned observation to generalize this result in the context of spin fill-ins satisfying the dominant energy condition (DEC). Following \cite{Bartnik5}, a $5$-uple $(\Sigma,\gamma,H,\alpha,\mathfrak{h})$ is called a spacetime Bartnik data set if $(\Sigma,\gamma)$ is an oriented, closed Riemannian manifold, $H$ and $\mathfrak{h}$ are smooth functions on $\Sigma$ and $\alpha$ is a smooth $1$-form on $\Sigma$. In this situation, a triplet $(\Omega,g,k)$ is a fill-in of such a spacetime Bartnik data if 
\begin{enumerate}
\item $(\Omega,g,k)$ is a compact initial data set, that is $(\Omega,g)$ is an $(n+1)$-dimensional compact Riemannian manifold with boundary and $k$ is a smooth symmetric $(0,2)$-tensor field on $\Omega$,

\item there exists an isometry $f:(\Sigma,\gamma)\rightarrow(\partial\Omega,g_{|\partial\Omega})$ such that
\begin{enumerate}
\item $f^\ast H_g=H$, where $H_g$ is the mean curvature of $\partial\Omega$ in $(\Omega,g)$ with respect to the outward unit normal $\widetilde{\nu}$,

\item $f^\ast\big(k(\widetilde{\nu},\cdot)^{T}\big)=\alpha$,  

\item $f^\ast\big(\Tr_{g_{|\partial\Omega}}\,k\big)=\mathfrak{h}$.
\end{enumerate}
\end{enumerate}
Here $\Tr_{g_{|\partial\Omega}}$ denotes the trace operator on $\partial\Omega$ and $\omega^T$ is the tangent part to $\Sigma$ of a $1$-form $\omega$ defined along $\Sigma$. In the following, we will omit the isometry $f$ in the identification between $\Sigma$ and $\partial\Omega$. Then, a fill-in $(\Omega,g,k)$ of $(\Sigma,\gamma,H,\alpha,\mathfrak{h})$ satisfies the dominant energy condition, or is a DEC fill-in, if 
\begin{eqnarray*}\label{DEC}
\mu\geq|J|_g
\end{eqnarray*}
where $\mu$ and $J$ are respectively the mass density and the current density defined by
\begin{eqnarray*}
\mu =  \frac{1}{2}\big(R_g+(\Tr_g\,k)^2-|k|^2_g\big) 
\end{eqnarray*}
and 
\begin{eqnarray*}
J = {\rm div}_g\big(k-(\Tr_g\,k)g\big).
\end{eqnarray*}
Here $R_g$ and $\Tr_g$ denote respectively the scalar curvature and the trace operator of $(\Omega,g)$. In this situation, a natural generalization of the mean curvature is given by the function
\begin{eqnarray*}
\mathcal{H}_g:=H_g-\sqrt{|k(\widetilde{\nu},\cdot)^{T}|_g^2+(\Tr_{g_{|\partial\Omega}}\,k)^2}
\end{eqnarray*}
which corresponds, for a spacetime Bartnik data, to 
\begin{eqnarray}\label{GeneralizerMeanCurvature}
\mathcal{H}:=f^\ast\mathcal{H}_g=H-\sqrt{|\alpha|_{\gamma}^2+\mathfrak{h}^2}.
\end{eqnarray} 
When $\Sigma$ is endowed with a spin structure, we will say that $(\Omega,g,k)$ is a DEC spin fill-in of the spacetime Bartnik data $(\Sigma,\gamma,H,\alpha,\mathfrak{h})$ if $(\Omega,g,k)$ is a DEC fill-in and if $\Omega$ is a spin manifold which induces the given spin structure on $\Sigma$. 
We then have the following result.
\begin{theorem}\label{MainThm}
Given any Riemannian metric $\gamma$ on an $n$-dimensional spin manifold $\Sigma$, there exists a constant $\mathcal{H}_0$, depending only on $\Sigma$, $\gamma$ and the spin structure of $\Sigma$, such that, if $\min_\Sigma \mathcal{H}\geq\mathcal{H}_0$, there do not exist DEC spin fill-ins of $(\Sigma,\gamma,H,\alpha,\mathfrak{h})$. 
\end{theorem}
This question was tackled by Tsang in \cite{Tsang} where several partial results are proved. We remark that if $k=0$ then $\alpha=0$ and $\mathfrak{h}=0$, and the DEC gives the nonnegativity of the scalar curvature so that a DEC fill-in of $(\Sigma,\gamma,H)$ corresponds to a NNSC fill-in. In a same way, if $k=cg$, $c\neq 0$, then $\alpha=0$ and $\mathfrak{h}=nc$, and the DEC condition implies that the scalar curvature of $(\Omega,g)$ is bounded from below by a negative constant, namely $R_g\geq-n(n+1)c^2$. These remarks imply that our result covers both of these cases. As mentioned above, the proof is a direct consequence of an eigenvalue estimate for the first eigenvalue of the Dirac operator of a spacetime Bartnik spin data which admits a DEC spin fill-in. This lower bound is stated as follow.   
\begin{theorem}\label{MainEstimate}
If a spacetime Bartnik spin data $(\Sigma,\gamma,H,\alpha,\mathfrak{h})$ admits a DEC spin fill-in $(\Omega,g,k)$ with $\mathcal{H}>0$, then the first eigenvalue $\lambda_1(\Sigma,\gamma)$ of the Dirac operator of $(\Sigma,\gamma)$ satisfies
\begin{eqnarray*}\label{BoundMeanCurvature}
\min_\Sigma \mathcal{H} \leq 2|\lambda_1(\Sigma,\gamma)|.
\end{eqnarray*}
\end{theorem}
One can consider $\mathcal{F}_{spin}(\Sigma,\gamma)$, the set of the DEC spin fill-ins of the Riemannian spin manifold $(\Sigma,\gamma)$ without specifying the data $H$, $\alpha$ and $\mathfrak{h}$. Then Theorem \ref{MainEstimate} implies that if $\mathcal{F}_{spin}(\Sigma,\gamma)\neq\emptyset$, it holds that
\begin{eqnarray*}
\sup_{(\Omega,g,k)\in\mathcal{F}_{spin}(\Sigma,\gamma)}\min_\Sigma\mathcal{H}\leq 2|\lambda_1(\Sigma,\gamma)|<\infty.
\end{eqnarray*}
The proof of Theorem \ref{MainEstimate} relies on spin geometry and we refer especially to \cite{BourguignonHijaziMilhoratMoroianu,Friedrich2,Ginoux} for more details on this subject. Let us briefly recall what we need here. Since $(\Omega,g)$ is a Riemannian spin manifold, there exists a smooth Hermitian vector bundle over $\Omega$, the spinor bundle, denoted by $\SB$, whose sections are called spinor fields. The Hermitian scalar product is denoted by $\<\,,\,\>$. Moreover, the tangent bundle $T\Omega$ acts on $\SB$ by Clifford multiplication $X\otimes \psi\mapsto c(X)\psi$ for any tangent vector fields $X$ and any spinor fields $\psi$. On the other hand, the Riemannian Levi-Civita connection $\nabla$ lifts to the so-called spin Levi-Civita connection, also denoted by $\nabla$, and defines a metric covariant derivative on $\SB$ that preserves the Clifford multiplication. A quadruplet $(\SB,c,\<\,,\,\>,\nabla)$ which satisfies the previous assumptions is usually referred to as a Dirac bundle. The Dirac operator is then the first order elliptic differential operator acting on $\SB$ defined by $D:=c\circ\nabla$. The spin structure on $\Omega$ induces, via a choice of an unit normal field to $\partial\Omega\simeq\Sigma$, a spin structure on $\Sigma$. This allows to define the {\it extrinsic} spinor bundle $\RSB:=\SB_{|\Sigma}$ over $\Sigma$ on which there exists a Clifford multiplication $\mult$ and a metric covariant derivative $\nb$. The quadruplet $(\RSB,\mult,\<\,,\,\>,\nb)$ is thus endowed with a Dirac bundle structure. Similarly, the extrinsic Dirac operator is defined by taking the Clifford trace of the covariant derivative $\nb$ that is $\D:=\mult\circ\nb$. It is by now well-known that this operator can be expressed using the Dirac operator $D_\Sigma$ of $(\Sigma,\gamma)$ endowed with the induced spin structure. What is important to us here is that the first nonnegative eigenvalue of the extrinsic Dirac operator $\D$ corresponds to $|\lambda_1(\Sigma,\gamma)|$, the absolute value of the first eigenvalue of $D_\Sigma$ and so it only depends on $(\Sigma,\gamma)$ and the spin structure on $\Sigma$. 

\vspace{0.2cm}

\noindent {\it Proof of Theorem \ref{MainEstimate}:}
Let $(\Omega,g,k)$ be a DEC spin fill-in of the spacetime Bartnik spin data $(\Sigma,\gamma,H,\alpha,\mathfrak{h})$ and consider the modified spin covariant derivatives defined by
\begin{eqnarray}\label{ModifiedConnection}
\nabla^\pm_X\psi:=\nabla_X\psi\pm\frac{i}{2}c\big(k(X)\big)\psi
\end{eqnarray}
for $X\in\Gamma(T\Omega)$ and $\psi\in\Gamma(\SB)$. The associated Dirac-type operators given by $D^\pm:=c\circ\nabla^\pm$ are easily seen to satisfy
\begin{eqnarray}\label{ModifiedDirac} 
D^\pm\psi=D\psi\mp\frac{i}{2}(\Tr_g\,k)\psi.
\end{eqnarray}
These are first order elliptic differential operators whose formal adjoints, with respect to the $L^2$ scalar product on $\SB$, are $D^\mp$ as deduced from the following integration by parts formulae
\begin{eqnarray}\label{IPP-Dirac}
\int_\Omega\<D^\pm\psi,\varphi\>d\mu=\int_\Omega\<\psi,D^\mp\varphi\>d\mu-\int_\Sigma\<c(\nu)\psi,\varphi\>d\sigma
\end{eqnarray} 
for all smooth spinor fields $\psi$, $\varphi$ on $\Omega$. Here $d\mu$ (resp. $d\sigma$) denotes the Riemannian volume (resp. area) element of $(\Omega,g)$ (resp. $(\Sigma,\gamma)$) and $\nu$ is the inner unit normal to $\Sigma$ in $(\Omega,g)$. In a same way, a straightforward computation implies that
\begin{eqnarray}\label{AdjointConnection}
\big(\nabla^\pm)^\ast\nabla^\pm\psi=-\sum_{j=1}^{n+1}\nabla^{\mp}_{e_j}\nabla^{\pm}_{e_j}\psi
\end{eqnarray}
where $\big(\nabla^\pm)^\ast$ denote the formal adjoints of the modified connection $\nabla^\pm$ and $\{e_1,\cdots,e_{n+1}\}$ is a local $g$-orthonormal frame of $T\Omega$. In particular, the Stokes formula leads to 
\begin{eqnarray}\label{IPP-Connection}
\int_\Omega\<\big(\nabla^\pm)^\ast\nabla^\pm\psi,\psi\>d\mu=\int_\Omega|\nabla^\pm\psi|^2d\mu+\int_\Sigma\<\nabla^\pm_\nu\psi,\psi\>d\sigma.
\end{eqnarray}
Then, it follows from the fact that $\big(D^\pm)^\ast=D^\mp$, from (\ref{AdjointConnection}) and the classical Schr\"odinger-Lichnerowicz formula  
\begin{eqnarray*}
D^2\psi=\nabla^\ast\nabla\psi+\frac{R_g}{4}\psi
\end{eqnarray*}
that 
\begin{eqnarray}\label{RiemannianWittenFormula}
\big(D^\pm\big)^\ast D^\pm\psi=\big(\nabla^\pm)^\ast\nabla^\pm\psi+\frac{1}{2}\big(\mu\psi\pm ic(J)\psi\big)
\end{eqnarray}
for all $\psi\in\Gamma(\SB)$. This is the $(n+1)$-dimensional Riemannian counterpart of the formula obtained by Witten \cite{Witten1} in his proof of the positive energy theorem. Now observe that since the point-wise symmetric endomorphism $J^\pm:=\pm ic(J)$ satisfies $\big(J^\pm\big)^2\psi=|J|^{2}_g\psi$ it holds that 
\begin{eqnarray*}
\<J^\pm\psi,\psi\>\geq-|J|_g|\psi|^2
\end{eqnarray*}   
so that the DEC and (\ref{RiemannianWittenFormula}) imply the point-wise inequalities
\begin{eqnarray*}
\<\big(D^\pm\big)^\ast D^\pm\psi,\psi\>\geq \<\big(\nabla^\pm)^\ast\nabla^\pm\psi,\psi\>
\end{eqnarray*}
for all $\psi\in\Gamma(\SB)$. Now integrating by part these estimates on $\Omega$ using (\ref{IPP-Dirac}) and (\ref{IPP-Connection}) leads to the following important integral inequalities
\begin{eqnarray}\label{IntegralInequality}
\int_\Omega\Big(|\nabla^\pm\psi|^2-|D^\pm\psi|^2\Big)d\mu\leq -\int_\Sigma\<\nabla^\pm_\nu\psi+c(\nu)D^\pm\psi,\psi\>d\sigma.
\end{eqnarray}
From the very definitions (\ref{ModifiedConnection}) and (\ref{ModifiedDirac}) of the modified covariant derivatives and the associated Dirac operators, we compute that 
\begin{eqnarray}\label{BoundaryExpression}
-\nabla^\pm_\nu\psi-c(\nu)D^\pm\psi =\D\psi-\frac{1}{2}\big(H\psi \mp ic(\mathcal{V})\psi\big)
\end{eqnarray}
where 
\begin{eqnarray*}
\mathcal{V}:=\alpha^\sharp+\mathfrak{h}\nu\in\Gamma(T\Omega_{|\Sigma})
\end{eqnarray*}
since $(\Omega,g,k)$ is a fill-in of the data $(\Sigma,\gamma,H,\alpha,\mathfrak{h})$. Here $\sharp:T^\ast \Omega\rightarrow T\Omega$ denotes the classical musical isomorphism between the cotangent bundle and the tangent bundle. Observe that the endomorphisms $\mathcal{V}^\pm:=\pm ic(\mathcal{V})$ of $\RSB$
is point-wise symmetric with respect to the Hermitian structure and satisfies $(\mathcal{V}^\pm)^2\psi=|\mathcal{V}|^2\psi$ in such a way that 
\begin{eqnarray}\label{MeanCurvatureBound}
\<\mathcal{V}^\pm\psi,\psi\>\geq-|\mathcal{V}|_g|\psi|^2=-\sqrt{\mathfrak{h}^2+|\alpha|_\gamma^2}\,|\psi|^2
\end{eqnarray}
for all $\psi\in\Gamma(\RSB)$. Combining (\ref{IntegralInequality}), (\ref{BoundaryExpression}) and (\ref{MeanCurvatureBound}) yields the following integral inequalities 
\begin{eqnarray}\label{IntegralInequalityFinal}
\int_\Omega\Big(|\nabla^\pm\psi|^2-|D^\pm\psi|^2\Big)d\mu\leq\int_\Sigma\<\D\psi-\frac{1}{2}\mathcal{H}\psi,\psi\>d\sigma.
\end{eqnarray}
which hold for all $\psi\in\Gamma(\SB)$ and where $\mathcal{H}$ is the generalized mean curvature function defined in (\ref{GeneralizerMeanCurvature}).

Now we are going to show that one can extend any spinor fields on $\Sigma$ in a suitable way. For this, we recall that the map $\chi:=ic(\nu)$ is a boundary chirality operator (in the sense of \cite[Example 7.26]{BarBallmann2}) and so it is an orthogonal involution of $\RSB$ which induces an orthogonal splitting $\RSB=\RSB^+\oplus\RSB^-$ into the eigenbundles of $\chi$ for the eigenvalues $\pm 1$. The associated projection maps $P^\pm:\RSB\rightarrow\RSB^\pm$ define elliptic local boundary conditions for the Dirac-type operators $D^\pm$ (see \cite[Corollary 7.23]{BarBallmann2} for example). This implies that the operators 
\begin{eqnarray*}
D^+_\pm:\big\{\psi\in H^1(\SB)\,/\,P^\pm\psi_{|\Sigma}=0\big\}\longrightarrow L^2(\SB)
\end{eqnarray*}
are of Fredholm type and that if $\Phi\in\Gamma(\SB)$ and $\Psi\in\Gamma(\RSB)$ are smooth spinor fields, any solutions $\psi\in\Gamma(\SB)$ of the boundary problem 
\begin{equation}\label{GeneralBVP}
\left\lbrace
\begin{array}{ll}
D^+\psi= \Phi & \text{ on }\Omega\\
P^\pm\psi_{|\Sigma} = P^\pm \Psi & \text{ on } \Sigma
\end{array}
\right.
\end{equation}
are smooth. The same holds for the operators $D^-_\pm$. It turns out that these operators are isomorphisms. To prove this fact, we notice that it is enough to show that $D^+_\pm$ and $D^-_\pm$ are one-to-one since it follows from the integration by parts formulae (\ref{IPP-Dirac}) that the adjoint of $D^+_\pm$ is $D^-_\mp$. So take $\psi_\pm$, non trivial, in the kernel of $D^+_\pm$, that is $\psi_\pm\in\Gamma(\SB)$ satisfies (\ref{GeneralBVP}) with $\Phi=0$ and $\Psi=0$. In particular, $\psi_\pm$ is smooth on $\Omega$. On the other hand, from the self-adjointness of the Dirac operator $\D$ and the fact that $\D\chi=-\chi\D$, we get that  
\begin{eqnarray}\label{DiracBoundary}
\int_\Sigma\<\D\psi,\psi\>d\sigma=2\int_\Sigma {\rm Re}\<\D(P^\pm\psi),P^\mp\psi\>d\sigma
\end{eqnarray}
for all $\psi\in\Gamma(\SB)$. Using this formula for $\psi=\psi_\pm$, we deduce from (\ref{IntegralInequalityFinal}) that 
\begin{eqnarray*}
0\leq\int_\Omega|\nabla^\pm\psi_\pm|^2d\mu\leq-\frac{1}{2}\min_\Sigma \mathcal{H}\int_\Sigma|\psi_\pm|^2d\sigma.
\end{eqnarray*}
Since we assumed that $\mathcal{H}$ is positive on $\Sigma$, we conclude that $\psi_\pm$ is zero on $\Sigma$ and $\nabla^\pm\psi_\pm=0$. This leads to a contradiction since this last property implies that $\psi_\pm$ is nowhere vanishing on $\Omega$. The same holds for $D^-_\pm$.

Now take $\Psi_1\in\Gamma(\RSB)$ be an eigenspinor for the operator $\D$ associated with the eigenvalue $|\lambda_1(\Sigma,\gamma)|$ and so, from the previous discussion, there exists an unique smooth solution $\varphi\in\Gamma(\SB)$ satisfying 
$$
\left\lbrace
\begin{array}{ll}
D^+\varphi= 0 & \text{ on }\Omega\\
P^+\varphi_{|\Sigma} = P^+ \Psi_1 & \text{ on } \Sigma.
\end{array}
\right.
$$
Taking $\psi=\varphi$ in (\ref{IntegralInequalityFinal}) and using the fact that
\begin{eqnarray*}
2{\rm Re}\<P^-\Psi_1,P^-\varphi\>\leq|P^-\Psi_1|^2+|P^-\varphi|^2,
\end{eqnarray*}
finally lead to 
\begin{eqnarray*}
0\leq \Big(|\lambda_1(\Sigma,\gamma)|-\frac{1}{2}\min_\Sigma\mathcal{H}\Big)\int_\Sigma |\varphi|^2d\sigma
\end{eqnarray*}
which implies the estimate of Theorem \ref{MainEstimate}. To get this last inequality, we implicitly used the fact that 
\begin{eqnarray*}
|\lambda_1(\Sigma,\gamma)|\int_\Sigma|P^-\Psi_1|^2d\sigma & = & \int_\Sigma\<\D(P^+\varphi),P^-\Psi_1\>d\sigma\\
& = & \int_\Sigma\<P^+\varphi,\D(P^-\Psi_1)\>d\sigma\\
& = & \int_\Sigma\<P^+\varphi,P^+(\D\Psi_1)\>d\sigma\\
& = & |\lambda_1(\Sigma,\gamma)|\int_\Sigma|P^+\varphi|^2d\sigma
\end{eqnarray*}
which follows from the identity (\ref{DiracBoundary}) and the self-adjointness of $\D$. 
\hfill$\square$

\vspace{0.2cm}

We conclude this note by noticing that this method is probably easily generalized to other situations (like the Einstein-Maxwell equations with nonnegative cosmological constant). 

%%%%%%%%%%%%%%%%%%%%%%%%%%%%%%%%%%%%%%%%%%%%%%%%%%%%%%%%%%%%%%%%%%%%%%%%%%%%%%%%%%%%

\bibliographystyle{alpha}   
\bibliography{BiblioDEC}     

%%%%%%%%%%%%%%%%%%%%%%%%%%%%%%%%%%%%%%%%%%%%%%%%%%%%%%%%%%%%%%%%%%%%%%%
\end{document}